\documentclass[12pt]{amsart}
\usepackage{amsmath,amsthm, amscd, amssymb, amsfonts}

\newcommand{\Z}{{\mathcal Z}}
\newcommand{\Zz}{{\mathbb Z}}

\newcommand{\R}{{\mathcal R}}
\newcommand{\C}{{\mathcal C}}

\newcommand{\Ss}{{\mathcal S}}

\newcommand\op{\operatorname{op}}
\newcommand\lcm{\operatorname{l.c.m. }}

\newcommand\vect{\operatorname{Vec}}
\newcommand\res{\operatorname{res}}
\newcommand\tr{\operatorname{tr}}
\newcommand\Rep{\operatorname{Rep}}

\newcommand\Opext{\operatorname{Opext}}

\newcommand\Aut{\operatorname{Aut}}

\newcommand{\fde}{{\triangleright}}
\newcommand{\fiz}{{\triangleleft}}

\newcommand\qexp{\operatorname{qexp}}

\numberwithin{equation}{section}\theoremstyle{plain}

\newtheorem{theorem}{Theorem}[section]
\newtheorem{lemma}[theorem]{Lemma}
\newtheorem{corollary}[theorem]{Corollary}
\newtheorem{conjecture}[theorem]{Conjecture}
\newtheorem{proposition}[theorem]{Proposition}

\theoremstyle{definition}

\newtheorem{example}[theorem]{Example}

\newtheorem{question}[theorem]{Question}

\theoremstyle{remark}
\newtheorem{remark}[theorem]{Remark}

\newcommand\id{\operatorname{id}}

\def\pf{\begin{proof}}
\def\epf{\end{proof}}

\theoremstyle{remark}

\begin{document}

\renewcommand{\baselinestretch}{1.2}
\renewcommand{\thefootnote}{}
\thispagestyle{empty}
\title[On the exponent of tensor categories]{On the exponent of
tensor categories coming from finite groups}
\author{Sonia Natale}
\address{Facultad de Matem\'atica, Astronom\'\i a y F\'\i sica
\newline \indent
Universidad Nacional de C\'ordoba
\newline
\indent CIEM -- CONICET
\newline
\indent (5000) Ciudad Universitaria
\newline
\indent
C\'ordoba, Argentina}
\email{natale@mate.uncor.edu \newline
\indent \emph{URL:}\/ http://www.mate.uncor.edu/natale}
\thanks{This work was partially supported by CONICET, Fundaci\' on
Antorchas, Agencia C\'ordoba Ciencia, ANPCyT    and Secyt (UNC)}
\date{November 4, 2005}

\begin{abstract} We describe the exponent of a group-theoretical fusion
category $\mathcal C = \mathcal C(G, \omega, F, \alpha)$
associated to a finite group $G$ in terms of group cohomology. We
show that the exponent of $\C$ divides both $e(\omega) \exp G$ and
$(\exp G)^2$, where $e(\omega)$ is the cohomological order of the
3-cocycle $\omega$. In particular $\exp \C$ divides $(\dim \C)^2$.
\end{abstract}

\maketitle

\section{Introduction and main results}

Throughout this note we shall work over an algebraically closed
base field $k$ of characteristic zero. The notion of
(quasi)exponent of a finite-dimensional Hopf algebra $H$ has been
introduced in a series of papers by Etingof and Gelaki
\cite{eg-exp, eg-qexp} extending previous work of Kashina
\cite{kashinayd, kashina}. By definition, the exponent of $H$ is
the least integer $N$ for which
$$m_N (\id \otimes \mathcal S^{-2} \otimes \dots \otimes \mathcal
S^{-2N + 2}) \Delta_N = \epsilon 1,$$ where $\Delta_N : H \to
H^{\otimes N}$ and $m_N: H^{\otimes N} \to H$ are the iterated
comultiplication and multiplication maps, respectively. This gives
a non-commutative analogue of the exponent of a group.

It was conjectured, in the context of semisimple Hopf algebras,
that the order of a certain power map divides the dimension of
$H$. In terms of the exponent, the conjecture can be stated as
follows:

\begin{conjecture}\label{conjetura} (\cite{kashinayd}.) If $H$ is a semisimple Hopf algebra over $k$,
then the exponent of $H$ divides the dimension of $H$. \end{conjecture}

This problem has an affirmative answer in a number of cases, but
the general answer is still not known. Etingof and Gelaki have
proved several basic and important properties and
characterizations of the exponent, in particular,  they have shown
that the exponent divides $(\dim H)^3$. One important property of
the exponent is its gauge invariance: that is, the exponent does
not depend on the Hopf algebra itself but only on its tensor
category of representations. Generalizing the definitions for
finite dimensional Hopf algebras, Etingof introduced the
quasi-exponent of a finite rigid tensor category $\C$ in
\cite{etingof}.

\medbreak The main goal of this paper is to describe the exponent
of a large class of semisimple Hopf algebras, which exhausts all
known examples, in terms of group cohomology. Actually, this class
consists not only of semisimple Hopf algebras but also of
semisimple quasi-Hopf algebras. Our results will imply that the
exponent of $H$ divides $(\dim H)^2$ for all $H$ in this class.

\medbreak Group-theoretical fusion categories where introduced by
Ostrik in \cite{ostrik}. Let $G$ be a finite group, and let
$\omega: G \times G \times G \to k^{\times}$ be a normalized
3-cocycle. Let also $F \subseteq G$ be a subgroup and $\tau: F
\times F \to k^{\times}$ a normalized 2-cochain
$\omega\vert_{F\times F \times F} = d\tau$. A
\emph{group-theoretical category} is a tensor category equivalent
to the category $\C(G, \omega, F, \alpha)$ of
$k_{\alpha}F$-bimodules in the category $\vect^G_{\omega}$ of
$G$-graded vector spaces with associativity given by $\omega$. A
(quasi)-Hopf algebra $H$ is called group theoretical if $\Rep H$
is.

\medbreak Recall that the \emph{global dimension} of $\C$, denoted
$\dim \C$, is the sum of squares of the categorical dimensions of
simple objects in $\C$. If $\Rep H \simeq \C(G, \omega, F,
\alpha)$, then $\dim \C = \dim H = |G|$.

It is an open question whether every semisimple Hopf algebra over
$k$ is group-theoretical or not \cite{ENO}. Every
group-theoretical category is equivalent to the
 category of representations of a quasi-Hopf algebra. The explicit
structure, up to gauge equivalence, of group theoretical
quasi-Hopf algebras was given in \cite{fs-indic}, where other
invariants, the Frobenius-Schur indicators, were computed in terms
of the group-theoretical data $G$, $\omega$, $F$, $\tau$.

Using a result of Schauenburg on the center of a bimodule
category,  it was shown in \cite{gp-ttic} that a quasi-Hopf
algebra $H$ is group theoretical if and only if its quantum double
is gauge equivalent to a Dijkgraaf-Pasquier-Roche quasi-Hopf
algebra $D^{\omega}G$ \cite{dpr}.

\medbreak In this paper we prove the following characterization of
the exponent of a group-theoretical category.

\begin{theorem}\label{car-exp} The exponent of $\C = \C(G, \omega, F, \alpha)$ divides the
modified exponent $$\exp_{\omega}G : = \text{l.c.m}
[e(\omega_g)|g|: \, g \in G],$$ and moreover $\exp \C =
\exp_{\omega}G$ in either of the following cases:
\begin{itemize} \item[(i)] $|G|$ is odd, \item[(ii)] $\C$ admits a fiber functor. \end{itemize}\end{theorem}

Here $\omega_g$ is the restriction of $\omega$ to the subgroup
generated by $g$ and $e(\omega_g)$ is its cohomological order.
Condition (ii) means that $\C$ is the category of representations
a Hopf algebra. Theorem \ref{car-exp} allows us to give necessary
and sufficient conditions for $\exp \C$ to divide $\dim \C$ when
the last is odd or $\C$ admits a fiber functor. See Theorems
\ref{ciclicos2}, \ref{nec-suf}.

\medbreak The following theorems are proved as a consequence of
this characterization.

\begin{theorem}\label{main2} The exponent of the twisted quantum
double $D^{\omega}G$ divides $(\exp G)^2$.
\end{theorem}

In particular, the exponent conjecture holds true for all
semisimple quasi-Hopf algebras which are gauge equivalent to a
twisted Drinfeld double $D^{\omega}G$. We also prove that for the
quasi-Hopf algebra $D^{\omega}G$, the order of the element $\beta$
divides the exponent of $G$.

\begin{theorem}\label{main} Let $\C \simeq \C(G,
\omega, F, \tau)$ be a group-theoretical fusion category. Then
\begin{itemize}\item[(i)] $\exp G$ divides $\exp \C$.
\item[(ii)]$\exp \C$ divides $e(\omega)\exp G$. In particular, $\exp \C$
divides $(\dim \C)^2$.
\item[(iii)] $\exp \C$ divides $(\exp G)^2$.
\item[(iv)] $\exp \C$ and $\dim \C$ have the same prime divisors.
\end{itemize} \end{theorem}

The proof relies on the characterization result in \cite{gp-ttic}.
We note that the statement in part (iv) has been recently
established in \cite{yyy} for any semisimple Hopf algebra $H$.

\medbreak In particular, it follows from Theorem \ref{main} that
the exponent of $A$ divides $(\dim A)^2$ for all bicrossed
products arising from exact factorizations of finite groups
\cite[Theorem 1.3]{gp-ttic} and all their twistings, that is, for
all semisimple Hopf algebras which are twist equivalent to some
$A$ that fits into an \emph{abelian} exact sequence
\begin{equation}\label{ab-ext}k \to k^{\Gamma} \to A \to kF \to
k,\end{equation} where $F$ and $\Gamma$ form a \emph{matched pair}
of finite groups; see \cite{ma-ext, ma-ext2}.

In this case, we show that $\exp A$ divides $\exp \Opext(F,
\Gamma)\exp G$, where $G = F \bowtie \Gamma$ is the factorizable
group determining the matched pair and $\Opext(F, \Gamma)$ is the
abelian group classifying all extensions \eqref{ab-ext}. See
Corollary \ref{main-cor}. Among semisimple Hopf algebras arising
from abelian extensions, the conjecture on the exponent had been
established under additional restrictions \cite{kashina}.

In the context of abelian extensions like \eqref{ab-ext} we also
obtain, as an application of results of Masuoka \cite{ma-ext2}, a
result that is of independent interest: we prove a Hopf algebra
generalization of the Schur-Zassenhauss Theorem for finite groups.
Namely, suppose that the orders of $\Gamma$ and $F$ are relatively
prime. Then, after twisting the multiplication and
comultiplication if necessary, $A$ is equivalent to the split
extension $k^{\Gamma}\# kF$. See Proposition
\ref{schur-zassenhauss}.

\medbreak The paper is organized as follows: in Sections
\ref{exponente} and \ref{double} we recall some properties of the
exponent of a fusion category and of quasi-Hopf algebras,
respectively. In Section \ref{dribbon} we review results of
Altschuler and Coste on Drinfeld and ribbon elements and prove,
under certain assumptions, some results on the powers of the
Drinfeld element that generalize those in \cite{eg-exp}. Finally,
in Section \ref{gt-fc}, we prove our main results, using the
ribbon element for a twisted Drinfeld double; Subsections
\ref{subs1} and \ref{subs2} concern particularly the Hopf algebra
case.

\medbreak {\bf Acknowledgements.} The author is grateful to S.
Montgomery for  interesting discussions on the exponent and its
properties. She also thanks J. Carlson for kind e-mail
correspondence.

\section{Exponent of a tensor category}\label{exponente}

In this section we recall the notion of (quasi-)exponent of a
finite tensor category introduced by Etingof in \cite[Section
6]{etingof}. This generalization is based on the results on the
(quasi-)exponent of a finite dimensional Hopf algebra found in the
papers \cite{eg-exp, eg-qexp}.

Let $\C$ be a finite rigid tensor category over $k$ and let
$\Z(\C)$ be its Drinfeld's center, which is a braided tensor
category with respect to a canonical braiding $\beta_{U V}: U
\otimes V \to V \otimes U$. In the paper \cite{etingof}, the
\emph{quasi-exponent} $\qexp \C$ of $\C$ is defined as the
smallest integer $N$ such that $(\beta^2)^N$ is unipotent in
$\Z(\C)$. Here, $\beta^2$ is the natural automorphism $\beta_{V
U}\beta_{U V}: U \otimes V \to U \otimes V$.

\medbreak The notion of quasi-exponent of a finite tensor category
gives, by restriction, a notion of quasi-exponent of a finite
dimensional quasi-Hopf algebra $H$ which by definition is a gauge
invariant of $H$.

\medbreak If $\C$ is a \emph{fusion} category, then the
quasi-exponent of $\C$ is called the \emph{exponent} of $\C$ and
denoted $\exp \C$. We shall be interested in fusion categories of
the form $\C = \Rep H$ where $H$ is a finite dimensional
semisimple quasi-Hopf algebra. In this case $\exp \C$ will be
called the \emph{exponent of $H$} and denoted $\exp H$.

It follows from \cite[Proposition 6.3]{etingof} that the exponent
of $\C = \Rep H$ satisfies the following:
\begin{flalign}\label{prop-1}& \exp \C = \exp \Z(\C); & \\
\label{prop-2}&\exp \C \, \text{equals the order of} \, \beta^2.&
\end{flalign}
Moreover,  \cite[Theorem 5.1]{etingof} implies that  $\exp \C$ is
finite.

\section{Semisimple quasi-Hopf algebras}\label{double}

Let $(H, \Delta, \epsilon, \phi, \mathcal S, \alpha, \beta)$ be a
finite dimensional semisimple quasi-Hopf algebra \cite{drinfeld}
(later on indicated by $H$ for short). Here, $\phi \in (H^{\otimes
3})^{\times}$ is the \emph{associator}, $\mathcal S: H \to
H^{\op}$ is the \emph{quasi-antipode} and $\alpha, \beta \in H$
are related to $\mathcal S$ by
\begin{equation}\label{ant-ab}\mathcal S(h_1) \alpha h_2 = \epsilon(h)\alpha, \qquad h_1 \beta \mathcal S(h_2) = \epsilon (h) \beta, \qquad \forall h\in H; \end{equation}
\begin{equation}\label{largas}\phi^{(1)} \beta \mathcal S(\phi^{(2)}) \alpha \phi^{(3)} = 1 = \mathcal S(\phi^{(-1)}) \alpha \phi^{(-2)} \beta \mathcal S(\phi^{(-3)}),
 \end{equation}
where $\phi = \phi^{(1)} \otimes \phi^{(2)} \otimes \phi^{(3)}$
and $\phi^{-1} = \phi^{(-1)} \otimes \phi^{(-2)} \otimes
\phi^{(-3)}$.

\medbreak The category $\Rep H = : \Rep (H, \phi)$ is a fusion
category of global dimension $\dim \C = \dim H$ with associativity
given by the action of $\phi$.

\medbreak Let $H_1$ and $H_2$ be finite dimensional semisimple
quasi-Hopf algebras. The categories $\Rep H_1$ and $\Rep H_2$ are
tensor equivalent if and only if $H_1$ and $H_2$ are gauge
equivalent \cite{et-gel}; that is, if and only if there exists a
 an invertible normalized
element  $F \in H_1 \otimes H_1$ (a \emph{gauge transformation})
such that $(H_1)_F$ and $H_2$ are isomorphic as quasi-bialgebras,
where $(H_1)_F$ is the quasi-Hopf algebra $(H_1, \Delta_F,
\epsilon, \phi_F, \mathcal S_F, \alpha_F, \beta_F)$, such that
\begin{equation*}\Delta_F(h) = F \Delta(h)F^{-1},  \qquad h \in H, \end{equation*}
\begin{equation*}\phi_F = (1 \otimes F) (\id \otimes \Delta)(F) \phi (\Delta \otimes \id)(F^{-1}) (F^{-1} \otimes 1),  \end{equation*}
\begin{equation*}\alpha_F = \mathcal S(F^{(-1)}) \alpha F^{(-2)},  \qquad \beta_F = F^{(1)} \beta \mathcal S(F^{(2)}); \end{equation*}
with the notation $F = F^{(1)}  \otimes F^{(2)}$, $F^{-1} = F^{(-1)} \otimes
F^{(-2)}$.

\medbreak There is also a notion of \emph{quasitriangular}
quasi-Hopf algebra, requiring the existence of an invertible
$R$-matrix $R\in H \otimes H$. When $H$ is quasitriangular the
category $\Rep H$ is a braided tensor category with braiding given
by the action of $R$.

\medbreak The center $\Z(\Rep H)$ is equivalent to the
representation category of the quantum double $D(H)$ \cite{majid,
hausser-nill}: this is a quasitriangular semisimple quasi-Hopf
algebra with underlying vector space $H^* \otimes H$ and canonical
R-matrix
\begin{equation}\label{rmatriz} \R = \sum_i h_i \otimes D(h^i),
\end{equation} where $(h_i)_i$ is a basis of $H$ and $(h^i)_i$ is
the dual basis.

\medbreak Since the element $\R_{21}\R \in D(H) \otimes D(H)$ implements the natural isomorphism $\beta^2$ in the category $\Rep D(H) \simeq \Z(\Rep H)$, the results in
\cite{etingof} imply the following lemma.

\begin{lemma}\label{ord-r2} The  order of $\R_{21}\R$ is finite and equals the exponent of $\Rep H$. \qed \end{lemma}

\subsection{Twisted quantum doubles}\label{ejemplo} Let $G$ be a finite group and
$\omega$ a normalized 3-cocycle on $G$. The identity element in
$G$ will be denoted by $e$. Let $H$ be the quasi-Hopf algebra
$(k^G, \omega)$ of $k$-valued functions on $G$ with associator
$\omega \in k^G \otimes k^G \otimes k^G$. Then the quantum double
of $H$ is a quasitriangular quasi-Hopf algebra isomorphic to the
Dijkgraaf-Pasquier-Roche quasi-Hopf algebra $D^{\omega}G$
\cite{majid, dpr}.

This quasi-Hopf algebra is defined on the vector space $k^G
\otimes kG$ as follows.  Consider the maps $\theta, \gamma: G
\times G \to (k^G)^{\times}$, $$\theta(x, y) = \sum_{g \in
G}\theta_g(x, y)e_g, \quad \gamma(x, y) = \sum_{g \in
G}\gamma_g(x, y)e_g,$$ where  $e_g \in k^G$ are the canonical
idempotents: $e_g(h) = \delta_{g, h}$, $g, h \in G$, and
\begin{equation}\label{formula-theta}\theta_g(x, y)  = \dfrac{\omega(g, x, y) \omega(x, y, (xy)^{-1}g (xy))}{\omega(x, x^{-1}gx, y)},\end{equation}
\begin{equation} \gamma_g(x, y)  = \dfrac{\omega(x, y, g) \omega(g, g^{-1}xg, g^{-1}yg)}{\omega(x, g, g^{-1}yg)}.  \end{equation}
Then $D^{\omega}G$ is as an algebra the crossed product
$k^G\#_{\theta}kG$, with respect to the adjoint action, and it is
the crossed product $k^G{}^{\gamma}\#kG$ as a coalgebra, with
respect to the trivial coaction.

A basis of $D^{\omega}G$ consists of the elements $e_g \# x$, $g,
x \in G$. The multiplication and comultiplication are explicitly
determined by
\begin{align}(e_g\# x)(e_h\# y) & = \delta_{g, xhx^{-1}} \theta_g(x, y) e_g \# xy, \\
\Delta(e_g\# x) & = \sum_{st = g} \gamma_x(s, t) \, e_s\# x \otimes e_t\# x. \end{align}
The unit element is $1 : = 1 \# e = \sum_{g \in G}e_g \# e$ and the counit and antipode are determined by
\begin{align}\epsilon(e_g\# x) & = \delta_{g, e},\\
\mathcal S(e_g\# x) & = \theta_{g^{-1}}(x, x^{-1})^{-1} \, \gamma_x(g, g^{-1})^{-1} \, e_{x^{-1}g^{-1}x}\# x^{-1},
\end{align}
with $\alpha = 1$ and $\beta = \sum_{g} \omega(g, g^{-1}, g) e_g$.

This is a quasitriangular quasi-Hopf algebra with associator and
R-matrix given, respectively, by
\begin{equation}\label{r-matrix}\phi = \sum_{a, b, c} \omega(a, b, c)^{-1} e_a \otimes e_b \otimes
e_c, \quad \R = \sum_g e_g \otimes g.\end{equation}

It is also known that in the twisted Drinfeld double, $\beta$ is
an invertible element with inverse $\mathcal S(\beta) =
\beta^{-1}$.  Moreover, $\beta$ implements $\mathcal S^2$ by
conjugation;  that is, for all elements $h \in D^{\omega}G$ we
have the relation
\begin{equation}\label{s2-tw}\mathcal S^2(h) = \beta^{-1} h \beta. \end{equation}

Recall that an element $a \in D^{\omega}G$ is  called
\emph{group-like} if $\Delta(a) = a \otimes a$. Axiom
\eqref{ant-ab} for the antipode, combined with the fact that
$\alpha = 1$ in $D^{\omega}G$, implies that the set of non-zero
group-like elements form a subgroup of the group of units of
$D^{\omega}G$, denoted $G(D^{\omega}G)$, and we have $\mathcal
S(a) = a^{-1}$, for all $a \in G(D^{\omega}G)$.

The group-like elements in twisted Drinfeld doubles have been
completely described in \cite[Proposition 3.2]{mason-ng}:  an
element $a \in D^{\omega}G$ is group-like if and only if there
exist elements $x \in G$ and $f \in k^G$ such that $$\gamma_x =
df, \quad \text{and} \quad a = f \# x.$$ Here $df: G \to k$
denotes the coboundary of $f$ given by $df(g, h) =
f(g)f(h)f(gh)^{-1}$.

\section{Drinfeld elements and ribbon elements}\label{dribbon}

Let $(H, \phi, R)$ be a quasitriangular quasi-Hopf algebra. Let $u \in H$ be the element defined by
\begin{align*}u : & = \Ss(\phi^{(-2)} \beta \Ss(\phi^{(-3)})) \, \Ss(R^j) \alpha R_j \phi^{(-1)} \\ & =
m_{21}(\id \otimes \mathcal S)\left((\alpha \otimes 1) R p\right),
\end{align*} where $p = p_R = \phi^{(-1)} \otimes \phi^{(-2)} \beta \mathcal
S(\phi^{(-3)}) \in H \otimes H$ is the special element related to
a canonical adjunction formula in $\Rep H$ \cite{drinfeld,
hausser-nill}.

The element $u$ has been introduced by Altschuler and Coste in
\cite{alt-coste} generalizing the Drinfeld element for
quasitriangular Hopf algebras \cite{drinfeld-qt}. It satisfies
$$\Ss^2(h) = uhu^{-1},$$ for all $h \in H$ \cite[Section
3]{alt-coste}. Here, and in what follows, we are using the
notation $R = R_j \otimes R^j$, assuming a summation symbol over
repeated indexes. The action of $u$ on finite dimensional
representations of $H$ gives a canonical isomorphism between the
dual and double dual functors.

\begin{remark}\label{3.9} The Drinfeld element $u$ satisfies the
equation \cite[(3.9)]{alt-coste}
$$\Ss(\alpha)u = \Ss(R^j) \alpha R_j.$$
In particular, when $\alpha = 1$, we get $u = \Ss(R^j)R_j$, which
coincides  with the formula given by Drinfeld in the Hopf algebra
case. \end{remark}

Let us denote  $\widetilde R = (\alpha \otimes 1) R p \in H
\otimes H$. Let $n \geq 1$ be an integer. Following \cite{eg-exp,
eg-qexp} we define an element $\widetilde R_n \in H \otimes H$ by
the formula
\begin{align*}\label{rn}\widetilde R_n & = \widetilde R (\id \otimes \Ss^2)(\widetilde R) \dots
(\id \otimes \Ss^{2n-2})(\widetilde R) \\
& = \widetilde R_{i_1}\widetilde R_{i_2} \dots \widetilde R_{i_n}
\otimes \widetilde R^{i_1}\Ss^2(\widetilde R^{i_2}) \dots
\Ss^{2n-2}(\widetilde R^{i_n}),
\end{align*}
where $\widetilde R = \widetilde R_{i_1}\otimes \widetilde R^{i_1}
= \widetilde R_{i_2}\otimes \widetilde R^{i_2} = \dots =
\widetilde R_{i_n}\otimes \widetilde R^{i_n}$.

\medbreak As in the Hopf algebra case, this element is related to
the powers of the Drinfeld element. The following lemma gives the
precise relation.

\begin{lemma}\label{r_n} We have $u^n = m_{21}(\id \otimes
\Ss)(\widetilde R_n)$. \end{lemma}

\begin{proof} The proof is by induction on $n$. If $n = 1$, there is nothing to prove.
Let $n \geq 2$. We have
\begin{align*}m_{21}(\id \otimes
\Ss)(\widetilde R_{n+1}) & = \Ss^{2n+1}(\widetilde R^{i_{n+1}})
\dots \Ss(\widetilde R^{i_1}) \,
\widetilde R_{i_1} \dots \widetilde R_{i_{n+1}} \\
& = \Ss^{2n+1}(\widetilde R^{i_{n+1}}) \dots \Ss^3(\widetilde
R^{i_2}) \, u \widetilde R_{i_2} \dots
\widetilde R_{i_{n+1}} \\
& = \Ss^{2n+1}(\widetilde R^{i_{n+1}}) \dots \Ss^3(\widetilde
R^{i_2}) \, \Ss^2(\widetilde R_{i_2})
\dots \Ss^2(\widetilde R_{i_{n+1}}) u \\
& = \Ss^2\left(\Ss^{2n-1}(\widetilde R^{i_{n+1}}) \dots
\Ss(\widetilde R^{i_2}) \,
\widetilde R_{i_2} \dots \widetilde R_{i_{n+1}}\right) u \\
& = \Ss^2(u^n)u = u^nu = u^{n+1},
\end{align*} by induction, and using that $\Ss^2(u) = u$.
This proves the lemma. \end{proof}

\begin{remark}\label{rn-particular} Suppose that  $\alpha = 1$. By Remark \ref{3.9}, $u = \Ss(R^j)R_j$.
Then it is not difficult to show by induction on $n$ that in this
case we have $u^n = m_{21}(\id \otimes \Ss)(R_n)$, where $R_n$ is
defined by
$$R_n = R (\id \otimes \Ss^2)(R) \dots (\id \otimes
\Ss^{2n-2})(R),$$ as in the Hopf algebra case. \end{remark}

Unlike in the semisimple Hopf algebra case, it may happen that $H$
is  semisimple but $u$ is not a ribbon element in $H$. Suppose
$\alpha$ is invertible in $H$. According to the definition in 4.1
and Remark on page 13 in \cite{alt-coste}, a ribbon element in $H$
is the same as a central element $v \in H$ satisfying
\begin{flalign}&v^2 = u\Ss(u);&\\
&\Ss(v) = v;&\\
&\epsilon(v) = 1;&\\
&\label{delta-v}\Delta(v) = (v \otimes v)(R_{21}R)^{-1} =
(R_{21}R)^{-1}(v \otimes v).&
\end{flalign}

\medbreak For a finite group $G$ with normalized 3-cocycle
$\omega$, the expression for the Drinfeld element in the twisted
quantum double $D^{\omega}G$ is the following:
\begin{equation}u = \sum_{g \in G}\omega(g, g^{-1}, g)^{-2} \, e_g \# g^{-1}. \end{equation}
Moreover, $D^{\omega}G$ is a ribbon quasi-Hopf algebra with ribbon element $v$ given by
\begin{equation}\label{formulav}v = \sum_{g \in G}\omega(g, g^{-1}, g)^{-1} \, e_g \# g^{-1}, \end{equation}
and the following relation holds:
\begin{equation}\label{u-v}v = \beta u.\end{equation} See \cite[Section 5]{alt-coste}.

The action of the ribbon element $v$ on irreducible
representations gives the (twisted) \emph{modular invariant}
matrix $T$ studied in various papers, see for instance
\cite{bantay, gcr, acm}.

\medbreak We also note the following simpler formula for the
inverse of the ribbon element:
\begin{equation}\label{inverse-v}v^{-1} = \sum_{g\in G} e_g \# g.
\end{equation}

\medbreak Recall the expression \eqref{r-matrix} for the canonical
$R$-matrix $\R \in D^{\omega}G \otimes D^{\omega}G$.

\begin{lemma}\label{formula} Let $n \geq 1$. Then $\R_n = \sum_{g \in G}
e_g \otimes (g\beta^{-1})^n\beta^{n}$. \end{lemma}

Note in addition the following expression for the $n$th power
$g^n$ in $D^{\omega}G$ of an element $g \in G$:
$$g^n = \sum_{s \in G} \theta_s(g, g)\theta_s(g^2, g) \dots \theta_s(g^{n-1}, g) \, e_s \# g^n.$$

\begin{proof} Using relation \eqref{s2-tw}, we compute
\begin{align*}\R_n & = \R_{i_1}\R_{i_2} \dots \R_{i_n} \otimes \R^{i_1}\Ss^2(\R^{i_2}) \dots
\Ss^{2n-2}(\R^{i_n}) \\
& =  \R_{i_1}\R_{i_2} \dots \R_{i_n} \otimes \R^{i_1} (\beta^{-1}
\R^{i_2} \beta)(\beta^{-2} \R^{i_2} \beta^2) \dots
\beta^{-n+1}\R^{i_n}\beta^{n-1} \\
& = \R_{i_1}\R_{i_2} \dots \R_{i_n} \otimes (\R^{i_1} \beta^{-1})
(\R^{i_2} \beta^{-1})  \dots (\R^{i_n}\beta^{-1})\beta^{n},
\end{align*} which, in view of \eqref{r-matrix}, equals
\begin{equation*}\sum_{g_1, \dots, g_n}e_{g_1}\dots e_{g_n} \otimes (g_1\beta^{-1}) \dots (g_n\beta^{-1}) \beta^n
= \sum_{g}e_{g} \otimes (g\beta^{-1})^n \beta^{n},
\end{equation*} as claimed. \end{proof}

\medbreak The exponent of $D^{\omega}G$ can be characterized in terms of the ribbon element \eqref{formulav}.

\begin{lemma}\label{caract-v} The exponent of $D^{\omega}G$ equals the smallest positive integer
$N$ such that $v^N \in G(D^{\omega}G)$.
\end{lemma}

\begin{proof} By Lemma \ref{ord-r2} the exponent of $D^{\omega}G$ equals the order of $\R_{21}\R$.
Since $v$  is a ribbon element for $D^{\omega}G$, the lemma follows in view of formula \eqref{delta-v}.
\end{proof}

\section{Exponent of group theoretical fusion categories}\label{gt-fc}

Let $G$ be a finite group, and let $F \subseteq G$ be a subgroup.
Let also  $\omega: G \times G \times G \to k^{\times}$ be a
normalized 3-cocycle, and $\tau: F \times F \to k^{\times}$ a
normalized 2-cochain, such that $\omega\vert_{F\times F \times F}
= d\tau$.

Consider the category $\vect^G_{\omega}$ of finite dimensional
$G$-graded vector spaces, with associativity constraint given by
$\omega$. That is, $\vect^G_{\omega}$ is the representation
category of the quasi-Hopf algebra $(k^G, \omega)$. Since the
twisted group algebra $k_{\tau}F$ is an algebra  in
$\vect^G_{\omega}$, the category $\C(G, \omega, F, \tau)$ of
$k_{\tau}F$-bimodules in $\vect^G_{\omega}$ is a tensor category.
A \emph{group theoretical} category is by definition a fusion
category equivalent to  $\C(G, \omega, F, \tau)$ for some $G$,
$F$, $\omega$, $\tau$ \cite[Section 3]{ostrik}.

\medbreak A (quasi)-Hopf
algebra $H$ is called group theoretical if $\Rep H$ is.
By the results in \cite{gp-ttic}, a quasi-Hopf algebra $H$ is
group theoretical if and only if its quantum double is gauge
equivalent to a Dijkgraaf-Pasquier-Roche quasi-Hopf algebra
$D^{\omega}G$. In particular, if $\Rep H \simeq \C(G, \omega, F,
\tau)$, then $\dim H = \vert G \vert$.

\medbreak Properties \eqref{prop-1} and \eqref{prop-2} of the exponent imply the following lemma.

\begin{lemma}\label{reduccion} Let $\C \simeq \C(G, \omega, F, \tau)$ be a group-theoretical category.
Then $\exp \C = \exp D^{\omega}G = \exp (k^G, \omega)$. \qed
\end{lemma}

In particular if $H$ is a group theoretical quasi-Hopf algebra,
then $\exp H = \exp D^{\omega}G = \exp (k^G, \omega)$, for
appropriate choice of a finite group $G$ and a normalized
3-cocycle $\omega$ on $G$, such that $\vert G \vert = \dim H$.

\begin{remark} Note that, by gauge invariance, the exponent of  $D^{\omega}G$
depends only  on the cohomology class of $\omega$. \end{remark}

It is well known that $|G|H^r(G, k^\times) = 0$, for all $r \geq
0$. However, this relation is \emph{not} always true if we replace
$|G|$ by $\exp G$. Nevertheless, we have the following weaker
annihilation property. This will be used next to prove a
divisibility property for the order of $\beta$ in a twisted
quantum double.

\begin{lemma}\label{normaliz} Let $N = \exp G$. There exists a normalized 3-cocycle $\widetilde \omega$
which is cohomologous to $\omega$ and such that
$$\widetilde\omega(g, g^{-1}, g)^N = 1, \quad \forall g \in G.$$
\end{lemma}

\begin{proof} Since $|\langle g \rangle|$ divides the exponent of $G$, for all $g \in G$, then  $\exp G$
anihilates $H^3(\langle g\rangle, k^{\times})$, for all $g \in G$.

Therefore, for all $g \in G$, there exists a normalized 2-cochain
$f^g: \langle g \rangle \times \langle g \rangle \to k^{\times}$
such that $$\omega(x, y, z)^N = df^g(x, y, z) = f^g(xy, z)f^g(x,
y) f^g(x, yz)^{-1}f^g(y, z)^{-1},$$ for all $x, y, z \in \langle g
\rangle$. Because $\langle g^{-1} \rangle = \langle g \rangle$, we
may choose $f^g$ in  a way such that $f^{g^{-1}} = f^g$.

Next we define a normalized 2-cochain $f: G \times G \to
k^{\times}$ in the form
$$f(g, h) = \begin{cases}f^g(g, h), \quad \text{if} \, h \in \langle g \rangle, \\ 1, \quad \text{otherwise}.
\end{cases}$$
Then, for all $g \in G$, we have
\begin{align*}(df)(g, g^{-1}, g) & = f(e, g) f(g, g^{-1}) f(g, e)^{-1} f(g^{-1}, g)^{-1}
\\ & = f(g, g^{-1}) f(g^{-1}, g)^{-1}
\\ & = f^g(g, g^{-1}) f^{g^{-1}}(g^{-1}, g)^{-1}  \\ & = df^g(g, g^{-1}, g)
\qquad (\text{since} \, f^g = f^{g^{-1}}) \\ & = \omega(g, g^{-1},
g)^N. \end{align*} The lemma is established by putting $\widetilde
\omega = \omega d(f^{-\frac{1}{N}})$. \end{proof}

As a consequence of Lemma \ref{normaliz} we obtain the following.

\begin{corollary} The  order of $\beta$ in the group of units of $D^{\omega}G$ divides the exponent of $G$.
\qed \end{corollary}

\medbreak Let $G$ be a finite group and let $\omega$ be a
normalized 3-cocycle  on $G$.  In what follows we shall give a
proof of the characterization in Theorem \ref{car-exp}.

\medbreak Let $n \geq 1$. We introduce a map $\pi_{n, \omega}: G
\to k^{\times}$, by the formula $\pi_1 = \epsilon$, and $\pi_{n,
\omega}(g) = \pi_{n-1, \omega}(g) \omega(g, g^{n-1}, g)$, $g \in
G$, $n \geq 2$. In other words,
\begin{equation}\label{def-pin}\pi_{n, \omega}(g) = \omega(g,
g^{n-1}, g) \omega(g, g^{n-2}, g) \dots \omega(g, g,
g).\end{equation} Compare with formula (A.3) in \cite{gcr}.

\medbreak Let $g \in G$ and suppose that $|g|$ divides $n$. The
following relation is easily seen and will be frequently used in
what follows:
\begin{equation}\label{power-pi} \pi_{n, \omega}(g) = \pi_{|g|,
\widetilde\omega}(g),
\end{equation}
where $\widetilde \omega: \langle g\rangle \times \langle g\rangle
\times \langle g\rangle \to k^{\times}$ is the 3-cocycle obtained
from $\omega^{\frac{n}{|g|}}$ by restriction.

\medbreak Recall the expression \eqref{inverse-v} for the inverse
of the ribbon element in $D^{\omega}G$. The following lemma
follows from a straightforward computation.

\begin{lemma}\label{potencia-v} Let $n \geq 1$. Then $v^{-n} =
\sum_{g \in G} \pi_n(g) \, e_g \# g^n$. \qed
\end{lemma}

Let us denote by $\widehat G$ the group of (one-dimensional)
characters on the group $G$. So that $\widehat G = G(k^G)$.

\begin{proposition}\label{basica} \begin{itemize}\item[(i)] The exponent of $G$ divides $\exp D^{\omega}G$.

\item[(ii)] The exponent of $D^{\omega}G$ equals the smallest integer $n$ with the properties
$$\pi_{n, \omega} \in \widehat G, \quad g^n = e, \quad \forall g \in
G.$$
\end{itemize}\end{proposition}

\begin{proof} By Lemma \ref{caract-v} the exponent of $D^{\omega}G$ equals the smallest positive integer $n$ such
that $v^n$ (hence also $v^{-n}$) belongs to the group
$G(D^{\omega}G)$.

Suppose that $v^{n} \in G(D^{\omega}G)$. By the description of
group-like elements in twisted Drinfeld doubles from
\cite[Proposition 3.2]{mason-ng} (c.f. Subsection \ref{ejemplo}),
there exist $x \in G$ and $f \in k^G$ such that $\gamma_x = df$
and
$$v^n = f \# x.$$ Applying the map $\epsilon \otimes \id$ to both
sides, we get $(\epsilon \otimes \id) (v^n) = \epsilon(f)x$. But
$(\epsilon \otimes \id) (v) = e$ because of formula
\eqref{formulav}, and $\epsilon \otimes \id: D^{\omega}G \to kG$
is an  algebra map. Therefore $\epsilon(f)x = e$. This implies
that $x = e$, whence $df = \gamma_e = 1$, and thus  $f \in
\widehat G$. In particular $v^n = f \# e \in k^G$ and $v^{-n} =
f^{-1} \# e$. It follows from Lemma \ref{potencia-v} that $g^n =
e$, for all $g \in G$. Therefore $\exp G$ divides $n$. Hence part
(i) follows.

\medbreak Using that $g^n = e$, for all $g \in G$, Lemma
\ref{potencia-v} gives
$$v^{-n} = \sum_{g \in G} \pi_n(g) e_g \# g^n = \left(\sum_{g \in G} \pi_n(g) e_g\right) \# e.$$
Therefore $\pi_n$ must be a character of $G$. This finishes the
proof of the Proposition.
\end{proof}

\begin{theorem}\label{conj-twdouble} The exponent of $D^{\omega}G$ divides $e(\omega) \exp
G$.  \end{theorem}

\begin{proof} We may assume that $\omega(x, y, z)^{e(\omega)} =
1$, for all $x, y, z \in G$. It follows from Equation
\eqref{power-pi} that
$$\pi_{e(\omega) \exp G, \omega} = 1.$$ This implies the theorem, in view of Proposition \ref{basica}.
\end{proof}

Since both $e(\omega)$ and $\exp G$ divide $|G|$, it follows from
Theorem \ref{conj-twdouble} that for a twisted Drinfeld double
$D^{\omega}G$ the exponent divides $|G|^2$ which equals the
dimension of $D^{\omega}G$. We shall see next that this is not
true in general for any quasi-Hopf algebra.

\begin{example}\label{ejemplo-ciclico} Let $G = C_n = \langle a : \, a^n = 1\rangle$ be a cyclic group
of \emph{odd} order $n$. The group $H^3(C_n, k^{\times})$ is
cyclic of order $n$ parametrized by the cohomology classes of the
3-cocycles $\omega = \omega_{\zeta}$ defined by
\begin{equation}\label{co-ciclico}\omega(a^i, a^j, a^l) = \zeta^{lq_{ij}},
\end{equation}
$0 \leq i, j, l \leq n-1$, where $\zeta \in k^{\times}$ are the
$n$th roots of 1, and $q_{ij} \in \Zz$ is the quotient of $i+j$ in
the division by $n$. Explicitly,
$$q_{ij} = \begin{cases}0, \quad \text{if} \; i+j \leq n-1, \\ 1, \quad \text{if} \; i+j \geq n,
\end{cases}$$ and
$$\omega(a^i, a^j, a^l) = \begin{cases}1, \quad \text{if} \; i+j \leq n-1, \\ \zeta^l, \quad \text{if} \; i+j \geq
n. \end{cases}$$ for all $0 \leq i, j, l \leq n-1$.

\begin{lemma}\label{exp-ciclico} Let $\omega$ be as in \eqref{co-ciclico}.
Then $\exp D^{\omega}C_n = |\zeta| n$.
\end{lemma}

In particular, the exponent of the quasi-Hopf algebra $(k^{G},
\omega)$ needs not divide the order of $G$ (= its dimension).

\begin{proof} Let $N = \exp D^{\omega}C_n$. By Theorem
\ref{conj-twdouble}, $N / e(\omega) \exp C_n = |\zeta| n$. On the
other hand, straightforward computations show that
\begin{align*}\pi_{n, \omega} (a) & = \prod_{j = 1}^{n-1} \omega(a, a^j, a) =
\zeta, \\ \pi_{n, \omega} (a^{n-1}) & = \prod_{i = 1}^{n-1}
\omega(a^{n-1}, a^{n-i}, a^{n-1}) = \zeta^{(n-1)^2} =
\zeta.\end{align*} Since $n/N$, we have
$$\pi_{N, \omega}(a) = \pi_{n, \omega^{\frac{N}{n}}} (a) =
\zeta^{\frac{N}{n}},$$ and $$\pi_{N, \omega}(a^{n-1}) = \pi_{n,
\omega^{\frac{N}{n}}} (a^{n-1}) = \zeta^{\frac{N}{n}}.$$ By
Proposition \ref{basica} we must have $\zeta^{\frac{N}{n}} =
\zeta^{-\frac{N}{n}}$, implying that $|\zeta|$ divides
$2\dfrac{N}{n}$. Since $|\zeta|$ is odd by assumption, we obtain
$|\zeta|n / N$. Hence $N = |\zeta|n$ as claimed. \end{proof}

In the case $n = 2$, however, the exponent of $(k\Zz_2, \omega)$,
where $\omega$ is the non-trivial cocycle given by $\omega(a, a,
a) = -1$, does equal the order of $\Zz_2$, because in this case
$\zeta  = \zeta^{-1}$. \end{example}

We now introduce a modified exponent $\exp_{\omega} G$ of a finite
group $G$ endowed with a 3-cocycle $\omega$ that will be useful to
describe the exponent of $D^{\omega}G$. Let us denote by
$\omega_g$ the restriction of $\omega$ to the subgroup generated
by $g \in G$. Let
\begin{equation}\label{exp-omega}\exp_{\omega} G = \lcm [e(\omega_g) |g|: \, g \in
G].\end{equation}

It is clear that $\exp G$ divides $\exp_{\omega}G$, and
$\exp_{\omega}G$ divides $e(\omega) \exp G$. For instance, when
$G$ is cyclic of order $n$, we have $\exp_{\omega}G = e(\omega)n$.

\begin{proposition}\label{ciclicos} The exponent of $D^{\omega}G$ divides $\exp_{\omega} G$.
Moreover, equality holds when $|G|$ is odd. \end{proposition}

\begin{proof} Let $n = \exp_{\omega} G$. Then $|g| / n$, for all $g \in G$, and
by \eqref{power-pi}, $$\pi_{n, \omega}(g) = \pi_{|g|, \widetilde
\omega},$$ for all $g \in G$, where $\widetilde \omega$ is the
restriction to $\langle g \rangle$ of the 3-cocycle
$\omega^{\frac{n}{|g|}}$. Since, by definition of
$\exp_{\omega}G$, $e(\omega_g)$ divides $\dfrac{n}{|g|}$, then
$\pi_{n, \omega}(g) = 1$, for all $g \in G$. By Proposition
\ref{basica}, $\exp D^{\omega}G$ divides $n$, as claimed.

\medbreak Suppose now that $|G|$ is odd, and let $N = \exp
D^{\omega}G$.  Let $g \in G$. As in the proof of Lemma
\ref{exp-ciclico} we find that the class of the restriction of
$\omega^{\frac{N}{|g|}}$ to the subgroup generated by $g$ is
trivial.

Therefore $e(\omega_g) |g| / N$, for all $g \in G$, implying that
$\exp_{\omega}G/ N$. This shows that $\exp D^{\omega}G =
\exp_{\omega}G$, when $|G|$ is odd, as claimed. \end{proof}

\subsection*{Proof of Theorem \ref{main2}}\label{exp-cuadrado}   If $C$ is a cyclic group,
then $|H^3(C, k^{\times})| = |C|$, c.f. Example
\ref{ejemplo-ciclico}. Then $e(\omega_g) / |g|$, for all $g \in
G$. Thus $\exp_{\omega}G = \lcm [e(\omega_g)|g|: \, g \in G]$
divides $(\exp G)^2$. This implies the theorem in view of
Proposition \ref{ciclicos}. \qed

\begin{example} Suppose that $|G| = p^n$ and $\exp G = p^k$, where $p$ is a prime number, and $2k \leq n$.
Then $\exp D^{\omega}G$ divides $|G|$.

In particular, if $G$ is an extraspecial $p$-group of order $|G|
> p^3$, then $\exp D^{\omega}G$ divides $|G|$.

\begin{proof} The first claim is clear from Corollary \ref{exp-cuadrado}. Suppose that $G$ is
extraspecial and $|G| > p^3$. Then $|G| = p^{1 + 2m}$, for some
integer $m \geq 2$, and by definition, there is a central
extension $0 \to \mathbb Z_p \to G \to (\mathbb Z_p)^{2m} \to 1$;
c.f. \cite[8.23]{asch}. Note that, in general, if $N$ is a normal
subgroup of $G$, then $\exp G$ divides $\exp N  \exp (G/N)$. Since
both the kernel and the quotient are of exponent $p$, this implies
that $\exp G$ divides $p^2$. Hence $(\exp G)^2 / p^4$ and this
divides $|G|$ because $|G|
> p^3$. \end{proof} \end{example}

\subsection*{Proof of Theorem \ref{main}} To prove the theorem we
shall combine the results in this section with Lemma
\ref{reduccion} that tells us that $\exp \C =  \exp D^{\omega}G$.

Part (i) follows from Proposition \ref{basica}, (ii) follows from
Theorem \ref{conj-twdouble}, (iii) follows from Theorem
\ref{main2}, and finally (iv) is a consequence of (i) and (iii).
\qed

\medbreak Let us denote by $[G: g]$ the index in $G$ of the
subgroup generated by $g \in G$.

\begin{theorem}\label{ciclicos2} Suppose that
\begin{itemize}\item[(i)]$e(\omega_g)$ divides $[G:
g]$, for all $g \in G$. \end{itemize} Then
\begin{itemize}\item[(ii)]$\exp D^{\omega}G$ divides$|G|$.\end{itemize}

If the order of $G$ is odd, then (i) is equivalent to (ii).
\end{theorem}

\begin{proof} The first claim follows immediately from
Proposition \ref{ciclicos}.

Now suppose that $|G|$ is odd and that Condition (ii) holds. As in
the proof of Lemma \ref{exp-ciclico} we find that the class of the
restriction of  $\omega^{[G: g]}$ to the subgroup generated by $g$
is trivial. Then $e(\omega_g)$ divides $[G: g]$, and (i) holds.
\end{proof}

\subsection{Group theoretical Hopf algebras}\label{subs1} A group theoretical category is the representation
category of a Hopf algebra if and only if it admits a fiber
functor.

Recall from \cite{ostrik} that fiber functors of the group
theoretical fusion category $\C(G, \omega, F, \alpha)$ are
classified by equivalence classes of subgroups $\Gamma \subseteq
G$ and 2-cocycles $\beta$ on $\Gamma$ such that
\begin{itemize}
\item[(1)] $\omega\vert_{\Gamma} = 1$;

\item[(2)] $G = F\Gamma$;

\item[(3)] the cocycle $\alpha\beta^{-1}$ is non-degenerate on $F\cap
\Gamma$. \end{itemize}

\bigbreak \emph{In what follows we shall assume that $\C \simeq
\C(G, \omega, F, \alpha)$ is a group-theoretical category
admitting a fiber functor.} That is, $\C \simeq \Rep A$, for some
group-theoretical Hopf algebra $A$.

\medbreak The exponent of a group theoretical Hopf algebra turns
out to have a simpler description in terms of twisted Drinfeld
double.

\begin{lemma}\label{pi-gt} Let $\pi_{n,
\omega}$ be as in \eqref{def-pin}, $n \geq 1$. The following
statements are equivalent:
\begin{itemize}\item[(i)] $\pi_{n, \omega} : G \to k^{\times}$ is
a group homomorphism;

\item[(ii)] $\pi_{n, \omega}(g) = 1$, for all $g \in G$.
\end{itemize} \end{lemma}

\begin{proof} We only need to show (i) $\Longrightarrow$ (ii).
Let $\Gamma \subseteq G$ be a subgroup giving rise to a fiber
functor. We may assume that $\omega\vert_{\Gamma} = 1$ and
$\omega\vert_{F} = 1$.

Let $g \in G$, so that $g$ writes in the form $g = x s$, $x \in
F$, $s \in \Gamma$. If $\pi_{n, \omega}$ is a group homomorphism,
then $\pi_{n, \omega}(g) = \pi_{n, \omega}(x)\pi_{n, \omega}(s) =
1$. \end{proof}

\begin{proposition} Let $A$ be a group theoretical Hopf algebra with $$\Rep A \simeq \C(G, \omega, F,
\alpha).$$ Then the exponent of $A$ equals the order of the ribbon
element $v$ in  $D^{\omega}G$. \qed
\end{proposition}

\begin{proof} By Lemma \ref{reduccion}, $\exp A = \exp
D^{\omega}G$. Combining Proposition \ref{basica} with Lemma
\ref{pi-gt}, we see that $\exp D^{\omega}G$ equals the smallest
integer $n$ such that $g^n = e$, for all $g \in G$, and $\pi_{n,
\omega} = 1$. This is exactly the order of $v$ in view of formula
\eqref{inverse-v}. \end{proof}

\begin{remark} Since the index of a subgroup annihilates the kernel of
the restriction map, we find that the following relation holds for
every group theoretical Hopf algebra:
\begin{equation}\label{indices}\exp D^{\omega}G / ([G: F]; [G:
\Gamma]) \exp G.\end{equation} \end{remark}

\medbreak We next prove that the results in Proposition
\ref{ciclicos} and Theorem \ref{ciclicos2} hold for
group-theoretical Hopf algebras, without the assumption that $G$
has odd order.

\begin{proposition}\label{prod-exp} We have $\exp
\C = \exp_{\omega} G$. \end{proposition}

\begin{proof} We know that $\exp \C = \exp D^{\omega}G$. The proof
is identical to the proof of Proposition \ref{ciclicos}, using
Lemma \ref{pi-gt}.  \end{proof}

\subsection*{Proof of Theorem \ref{car-exp}} By Lemma \ref{reduccion}, $\exp \C = \exp
D^{\omega}G$. The theorem follows from Proposition \ref{ciclicos}
in case (i), and from Proposition \ref{prod-exp} in case (ii).
\qed

\begin{remark} Note that, in general, it is not true that $\exp \C$ divides
$\exp G$. For instance, if $\C = \Rep H_8$, where $H_8$ is the
8-dimensional Kac-Paljutkin Hopf algebra \cite{k-p}, then $G =
D_4$ the dihedral group of order $8$. We have in this case $\exp
\C = 8$ while $\exp G = 4$.\end{remark}

We next give some necessary and sufficient conditions for $\exp
\C$ to divide $\dim \C$.

\begin{theorem}\label{nec-suf} The following are equivalent:
\begin{itemize}\item[(i)] $\exp \C$ divides $\dim \C$;

\item[(ii)] $e(\omega_g)$
divides $[G: g]$, for all $g \in G$.
\end{itemize} \end{theorem}

\begin{proof} Identical to the proof of Theorem \ref{ciclicos2}, using Lemma \ref{pi-gt} and the
fact that $\exp \C = \exp D^{\omega}G$. \end{proof}

\begin{lemma} Suppose that $\exp
\C$ divides $\dim \C$. Then $e(\omega)$ divides $[G: g]^2$, for
all $g \in G$. \end{lemma}

\begin{proof}
By Lemma \ref{pi-gt}, $\pi_{|G|, \omega}(g) = 1$, for all $g \in
G$. Using Equation \eqref{power-pi},  we have $1 = \pi_{|G|,
\omega}(g) = \pi_{e(g), \widetilde\omega}(g)$, where $\widetilde
\omega: \langle g\rangle \times \langle g\rangle \times \langle
g\rangle \to k^{\times}$ is the 3-cocycle obtained from
$\omega^{[G: g]}$ by restriction.

As in the proof of Lemma \ref{exp-ciclico} we find that the class
of the restriction of $\omega^{[G: g]}$ to the subgroup generated
by $g$ is trivial.

The composition $$H^3(G, k^{\times}) \overset{\res}\to H^3(\langle
g \rangle, k^{\times}) \overset{\tr}\to H^3(G, k^{\times}),$$
where $\tr$ denotes the transfer map, is multiplication by the
index $[G: g]$ \cite{CE}. Therefore we find that $\omega^{[G:
g]^2} = 1$. Hence $e(\omega) / [G: g]^2$,  as claimed. \end{proof}

\medbreak Suppose that $G = F\Gamma$ is any factorizable finite
group. Let ${\widetilde H}^3(G, k^{\times})$ denote the kernel of
the restriction map $H^3(G, k^{\times}) \to H^3(F, k^{\times})
\oplus H^3(\Gamma, k^{\times})$. The following question is of a
purely cohomological nature.

\begin{question} Does the product $\exp {\widetilde H}^3(G, k^{\times}) \, \exp
G$ divide the order of $G$? \end{question}

An affirmative answer to this question would guarantee that the
exponent conjecture holds true for all group-theoretical Hopf
algebras.

\subsection{Abelian extensions}\label{subs2} The class of group-theoretical quasi-Hopf algebras
contains in particular the class of abelian bicrossed product Hopf
algebras, first studied by G. I Kac \cite{kac}. We refer the
reader to \cite{ma-ext, ma-ext2} for the main features of the
subject.

In what follows we shall consider a  fixed matched pair of finite
groups $(F, \Gamma)$ with respect to compatible actions $\fde :
\Gamma \times F \to F$, $\fiz : \Gamma \times F \to \Gamma$. These
actions determine a unique group structure on the product of $F$
with $\Gamma$, denoted $G = F \bowtie \Gamma$, in such a way that
$G$ admits an \emph{exact factorization} $G = F\Gamma$, $F \cap
\Gamma = 1$.

\begin{remark}\label{suc-esp} For every group with an exact factorization as
above, there are two convergent spectral sequences
$$H^p(F, H^q(\Gamma, k^{\times})) \Rightarrow H^{p+q}(G,
k^{\times}),$$ $$H^p(\Gamma, H^q(F, k^{\times})) \Rightarrow
H^{p+q}(G, k^{\times}).$$ \end{remark} These spectral sequences
come from the double complex in \cite{ma-ext, ma-ext2} whose total
complex gives a free resolution of the $G$-module $\mathbb Z$.

\medbreak For every class $(\sigma, \tau)$ in $\Opext (k^G, kF)$;
that is, $\sigma: F \times F \to (k^\Gamma)^{\times}$ and $\tau:
\Gamma \times \Gamma  \to (k^F)^{\times}$ are normalized
2-cocycles subject to certain compatibility conditions, there is a
bicrossed product Hopf algebra $A : = k^{\Gamma} \,
{}^{\tau}\#_{\sigma}kF$. This gives a one-to-one correspondence
between the equivalence classes of Hopf algebra extensions
\begin{equation}\label{abeliana}k \to k^\Gamma \to A \to kF \to k,\end{equation} affording the actions
$\fde$, $\fiz$, and the abelian group $\Opext (k^{\Gamma}, kF)$.

The \emph{Kac exact sequence} associated to the matched pair $(F,
\Gamma)$ \cite{kac, ma-ext} has the following form:
\begin{align*}
0 & \to H^1(G, k^{\times}) \xrightarrow{\res}   H^1(F, k^{\times})
\oplus  H^1(\Gamma, k^{\times}) \to \Aut(k^\Gamma \# kF) \\
 & \to  H^2(G, k^{\times}) \xrightarrow{\res}  H^2(F, k^{\times}) \oplus  H^2(\Gamma, k^{\times})
 \xrightarrow{\delta} \Opext(k^{\Gamma}, kF)  \\ & \to  H^3(G, k^{\times})
\xrightarrow{\res}  H^3(F, k^{\times}) \oplus H^3(\Gamma,
k^{\times}) \to \dots
\end{align*} where $\res$ denote the restriction maps.

\medbreak Hopf algebras $A$ arising from abelian exact sequences
are always group-theoretical: indeed,  $\Rep A \simeq \C(G,
\omega, F, 1)$, where $G = F \bowtie \Gamma$ and  $\omega$ is the
image of $(\sigma, \tau)$ in Kac exact sequence; see \cite[Theorem
1.3]{gp-ttic}. Explicitly, the 3-cocycle $\omega = \omega(\sigma,
\tau)$ can be represented by
\begin{equation}\omega (xg, x'g', x''g'') = \tau_{x''}(g \fiz x', g') \, \sigma_{g}(x', g' \fde x''), \end{equation}
for all $x, x', x'' \in F$, $g, g', g'' \in \Gamma$.

\medbreak Recall the Schur-Zassenhauss Theorem for finite groups
that says that any extension $G$ of a group $F$ by a group
$\Gamma$, with $|F|$ and $|\Gamma|$ relatively prime, splits. That
is, $G$ is a semidirect product of $F \rtimes \Gamma$. The
following proposition gives an analogue of this result for Hopf
algebras.

\begin{proposition}\label{schur-zassenhauss} Suppose that $|F|$ and $|\Gamma|$ are
relatively prime. Let $A$ be a Hopf algebra fitting into an
extension \eqref{abeliana}. Then $A$ is obtained from the split
extension $k^{\Gamma}\# kF$ by twisting the multiplication and the
comultiplication.
\end{proposition}

\begin{proof} Let $(\sigma, \tau)$ be the element in $\Opext (F, \Gamma)$ corresponding to the
extension \eqref{abeliana}. The indexes $[G, \Gamma] = |F|$ and
$[G, F] = |\Gamma|$ annihilate the kernel of the restriction map
$$H^3(G, k^{\times}) \to H^3(\Gamma, k^{\times}) \oplus H^3(F,
k^{\times}),$$ whence by exactness of the Kac sequence, $(\sigma,
\tau)$ belongs to the image of $\delta$. Now the result of Masuoka
on cocycle twists of bicrossed products \cite{ma-ext2} implies the
proposition. \end{proof}

\begin{corollary}\label{sz-cor}  Let the exact sequence \eqref{abeliana} and suppose that $|F|$, $|\Gamma|$ are relatively
prime.  Then $\exp A = \exp F \bowtie \Gamma$. \qed
\end{corollary}

In particular $A$ and all their cocycle twists satisfy the
exponent conjecture \ref{conjetura}.

\begin{proof} By Proposition \ref{sz-cor} and the twist invariance of the exponent, we may
assume that the extension \eqref{abeliana} splits. In this case
the 3-cocycle $\omega$ associated to $A$ under the Kac exact
sequence is trivial. The corollary follows from Theorem
\ref{main}. \end{proof}

The following is a consequence of Theorem \ref{main}.

\begin{corollary}\label{main-cor} Let $A$ be a Hopf algebra which fits into an abelian exact sequence
\eqref{abeliana}.  Then $\exp A$ divides $\exp \Opext(\Gamma, F)
\exp G$. \end{corollary}

\begin{proof} By Theorem \ref{main} (ii), $\exp A$ divides $e(\omega) \exp
G$, where $\omega$ is the 3-cocycle coming from the element in
$\Opext(F, \Gamma)$ corresponding to $A$ under the map $\delta$ in
the Kac sequence. The corollary follows from the exactness of the
sequence. \end{proof}

\medbreak The following proposition is a refinement of the
relation \eqref{indices} in the case of abelian exact sequences.
It generalizes the statement in Corollary \ref{schur-zassenhauss}.

\begin{proposition}Let $A$ be a Hopf algebra which fits into an abelian exact sequence
\eqref{abeliana}.  Then $\exp A$ divides $(|F|; |\Gamma|) \, \exp
G$. \qed \end{proposition}

\begin{proof} The indexes $[G, \Gamma] =
|F|$ and $[G, F] = |\Gamma|$ annihilate the kernel of the
restriction map $H^3(G, k^{\times}) \to H^3(\Gamma, k^{\times})
\oplus H^3(F, k^{\times})$. Hence $e(\omega)/(|F|; |\Gamma|)$.
This implies the proposition.
\end{proof}

\end{document}